\documentclass[10pt,reqno]{siamltex1213}
\usepackage{bm}
\usepackage[numbers]{natbib}
\usepackage{graphics,graphicx,color}
\usepackage{latexsym, amssymb, enumerate, amsmath,caption,subcaption}

\usepackage{url}
\usepackage{color}

\newcommand{\eps}{{\varepsilon}}

\renewcommand*{\Xi}{{\boldsymbol{\xi}}}

\renewcommand*{\O}{{\mathcal{O}}}

\newcommand{\dt}[0]{{\delta t}}
\newcommand{\Dt}[0]{{\Delta t}}

\newcommand{\tDt}[0]{\,\widetilde{\!\Delta \,}{\! t}}

\allowdisplaybreaks

\title{A note on implementations of the Boosting Algorithm and Heterogeneous Multiscale Methods}

\author{John Maclean\thanks{School of Mathematics and Statistics,  
University of Sydney, 
NSW 2006,
Australia (\email{J.Maclean@maths.usyd.edu.au)}}}

\begin{document}
\pagestyle{myheadings}\thispagestyle{plain} \markboth{A note on implementations of the BA and HMM}{John Maclean}

\maketitle

\begin{abstract}
We present improved convergence results for the Boosting Algorithm (BA), and demonstrate that an existing formulation of the Heterogeneous Multiscale Methods (HMM) is accurate to first order only in the macro time step, regardless of the order of the numerical solvers employed. These results are obtained by considering the BA and two other formulations of HMM as special cases of a general formulation of HMM applied to dissipative stiff ordinary differential equations. 
 \end{abstract}



\section*{Introduction} 
Many problems in the natural sciences are modelled by time-scale separated differential equations in which the dynamics of interest takes place on a slow time scale. The Heterogeneous Multiscale Methods (HMM) introduced in \cite{EEngquist03} provide the conceptual framework of micro and macro solvers to resolve the behaviour of the dynamics of interest without requiring a small integration time step. The micro solver performs short, fine-scale computations in the fast variables with micro time steps. The output from the micro solver is used to propagate the slow variables for macro time steps in the macro solver. 
\\ 
{These methods have been applied to dissipative and to highly oscillatory problems \cite{E03,ArielEtAl}, as well as multi-scale stochastic problems and partial differential equations \cite{VandenEijnden03,Abdulle09}, and to problems exhibiting spatial-scale separation; see \cite{EEtAl07} for a review.}
\\

\noindent
There are some differences in the literature on the exact manner in which the HMM macro and micro solver interact. In the formulation \cite{EEtAl07, VandenEijnden03, VandenEijnden07} the micro solver is applied during the macro solver whenever a vector field evaluation is required; for example, each time step of a fourth order Runge-Kutta macro solver would require four applications of the micro solver. We denote this approach by HMM1.
\\
 {A different approach to HMM is described in the more recent papers \cite{EEtAl09,E11,AbdulleEtAl12}, and employed in, for instance, \cite{LockerbyEtAl13}. HMM in these papers uses the micro and macro solver sequentially; for example, each time step of a fourth order Runge-Kutta macro solver would be preceded by an application of the micro solver, and would involve no additional applications of the micro solver. We denote this approach by HMM2.}\\
  {We remark that one cannot distinguish between HMM1 and HMM2 in some of the literature on HMM, as the two formulations are identical for the most common choice of macro solver, the forward Euler method. In this paper we will implement and analyse HMM1 and HMM2 with a Runge-Kutta macro solver in order to establish the important difference between them. \\}

\noindent
In practice, the large macro time steps of HMM may require that the fast variables be re-initialised in a manner consistent with the slow variables after each application of the macro solver. To mitigate this problem, the Boosting Algorithm (BA) was proposed in \cite{EEtAl09} {as an extension of the similarly motivated scheme in \cite{Ren07}.} {The BA alternates between single applications of the micro and macro solvers, with a smaller macro step than HMM. {This method was developed in \cite{Ren07,EEtAl09} for multi-scale ordinary, stochastic and partial differential equations, and extended to highly oscillatory systems in \cite{LeeEngquist14}.} In the context of dissipative multi-scale systems, the underlying assumption of the BA is that the slow manifold does not change too much over the smaller macro steps, so that the fast variables do not have to be expensively re-initialised or re-equilibrated by longer applications of the micro solver as done in HMM.}\\


\noindent
{A related class of methods called projective methods, proposed in \cite{GearKevrekidis03}, {has also been applied to time-scale separated ordinary, stochastic and partial differential equations, and to bifurcation analysis \cite{KevrekidisSamaey09,GivonEtAl06}. }{Formulations of these methods exist that resemble HMM2 \cite{GearKevrekidis03, LustRooseVandekerckhove06, RooseVandekerckhove06} and HMM1 \cite{GearLee05,LafitteEtAl13,MacleanGottwald15}. Unlike HMM, in projective methods the micro and macro solvers typically propagate both fast and slow variables. This allows projective methods to numerically integrate systems in which the slow and fast variables are not known; f}or a discussion of the benefits of this approach, see \cite{VandenEijnden07}.\\


\noindent
{There is an existing body of work on the convergence of HMM and the BA for ordinary differential equations. Of relevance are \cite{E03}, in which convergence results are given for HMM with an arbitrary order macro solver, albeit without proof or explicit formulation of the numerical scheme, and \cite{EEtAl09}, in which convergence results are given for the BA.}\\
{In this paper we provide rigorous convergence results for HMM1, HMM2 and the BA applied to stiff dissipative systems. We confirm for the HMM1 formulation the bound stated in \cite{E03}, and show that the HMM2 formulation incurs linear error in the macro step regardless of the order of the macro solver. The bound we derive for the BA is tighter than the existing bound in \cite{EEtAl09}; in particular, our bound removes a term that implies a restrictive condition should be satisfied by the micro solver.}\\

\noindent
The paper is organised as follows. In Section~\ref{sec:model} we consider a simple dissipative system. In Section~\ref{sec:methods} we present a general formulation of HMM, with Runge-Kutta macro and micro solvers and a variable number of micro steps per application of the micro solver. We then consider the BA, HMM1 and HMM2 as special cases of this general formulation. In Section~\ref{sec:analysis} we establish the convergence of the general formulation of HMM to the continuous solution of the slow processes in a multiscale system of deterministic ordinary differential equations by a simple extension of the convergence proof in \cite{E03}. We then apply this bound to the BA, HMM1 and HMM2. In Section~\ref{sec:numerics} we present results from numerical simulations corroborating our analytical findings. We conclude with a discussion in Section~\ref{sec:discussion}.

\section{Model}
\label{sec:model}
\noindent
For simplicity of exposition, we restrict our attention to systems with one slow variable $x_\eps$ and one fast variable $y_\eps$. We consider the deterministic multiscale system 
\begin{align}
\dot x_\eps &=f(x_\eps,y_\eps)\; ,
\label{baseslow}
\\
\dot y_\eps &=\frac{1}{\varepsilon}(-  \, y_\eps+ h_0(x_\eps))\; .
\label{basefast}
\end{align}
{We assume that there exists a slow manifold $x=h_\eps(y) = h_0(y) + \mathcal{O}(\eps)$, towards which initial conditions are attracted exponentially fast. On the slow manifold, the dynamics slows down and is approximately determined by}
\begin{align}
\dot X = F(X)\; ,
\label{e.CMT}
\end{align}
with reduced dynamics $X=x_\eps+{\mathcal{O}}(\eps)$. The reduced slow vector field is given by 
\begin{align}
\nonumber F(x) &= f(x,h_\eps(x)) \;  \\
\label{e.G0} &= f(x,h_0(x)) + \mathcal{O}(\eps) \; ,
\end{align}
where we identify the zeroth order approximation of the reduced vector field in \eqref{e.G0} because, as we shall see, the HMM and BA approximate $f(x,h_0(x))$ rather than $F(x)$. \\

\section{Methods}\label{sec:methods}
\noindent
We consider three numerical multiscale methods. Each method intermittently resolves the fast dynamics with a micro solver for micro time steps $\dt$ to relax the fast variable to the slow manifold, enabling the macro solver to resolve the slow dynamics with larger time steps. We select as the micro solver an explicit numerical method of order p. The three methods are conceptually described as follows:


\begin{remunerate}
\item \textbf{BA:} The macro solver is a numerical method of order P with macro step $\tDt \gg \dt $. In this approach the micro solver is applied for a single micro step before each time step of the macro solver. {As explained in \cite{EEtAl09}, one way to think of the BA is as a rescaling of time. In particular, it is equivalent to replacing $\eps$ in \eqref{basefast} with $\tDt\eps/\dt$ and then solving the new system of ordinary differential equations by standard numerical methods with time step $\tDt$. Since $\tDt\eps/\dt \gg \eps$, the new system is much less stiff and the BA presents a large computational advantage over solving \eqref{baseslow}--\eqref{basefast} directly.}   
\item \textbf{HMM1:} The macro solver is constructed by modifying a numerical method of order P to include an application of the micro solver before every vector field evaluation of \eqref{baseslow}. The macro step is $\Dt \gg \tDt$, and the micro solver is applied for sufficiently many steps to relax the fast variable.  
\item \textbf{HMM2:} The macro solver is given by a numerical method of order P with macro step $\Dt \gg \tDt$. The micro solver is applied for sufficiently many steps to relax the fast variable between each application of the macro solver.  
\end{remunerate}

\noindent
We will compare the three methods, focussing on a detailed analysis of the BA and the differences between HMM1 and HMM2. To do so, we consider a general formulation of HMM which allows for the construction of the BA, HMM1 and HMM2 formulations with particular parameter choices.

\noindent 
Throughout the paper, we use superscripts to denote elements of the macro solver and subscripts to denote elements of the micro solver. 
\subsection{General formulation of HMM}
We describe an HMM scheme in which the macro solver is given by a Runge-Kutta method of order P, except that the micro solver is applied for a variable number of steps before each vector field evaluation of the slow dynamics. Denote by $\varphi_{m,\dt}$ the flow map for the micro solver applied to the fast dynamics \eqref{basefast} with fixed slow variable for $m$ micro steps with time step $\dt$. Denote by $x^n$ and $y^n$ the approximations to the slow and fast variables at macro time $t^n$. The macro solver is then given by a weighted sum of increments, each of which cover a time step $\Dt$. The increments are generated recursively by
\begin{align}
\label{kj} {k}^{(j)}(x^{n},y^n)& = \Dt \, f (x^n + a^{(j)} k^{(j-1)},\, y^{n,j}_{M_j} ) \;\; ,
\end{align}
for $j=1,2,\dots,P$, where we define $y^{n,j}_{m}$ for $j=1,2,\dots,P$, $m=1,2,\dots,M_j$, as the output of the micro solver
\begin{align}
\label{ynjm} y^{n,j}_{m} &= \varphi_{m,\dt}\left(x^n + a^{(j)} k^{(j-1)},\, y^{n,j}_{0}\right) \;\; ,
\end{align}
with initial condition at $j=1$
\begin{align}
y^{n,1}_{0} &= y^n \;\;,
\end{align}
and at $j>1$
\begin{align}
\label{yj0} y^{n,j}_{0} &= y^{n,1}_{M_1} \;\;.
\end{align}
The nodes $a^{(j)}$ in \eqref{kj}, \eqref{ynjm} are those used in the increments of a Runge-Kutta solver of order P; in particular $0\le a^{(j)} \le 1$, and $a^{(1)} = 0$ always so that $k^{(1)}$ is defined explicitly. \\
\noindent The macro solver is then given by the weighted sum
\begin{align}
\label{macrox} x^{n+1} &= x^{n} + \sum_{j=1}^P b^{(j)} {k}^{(j)}(x^{n},y^{n}) \;\; ,
\end{align}
%
where the weights $b^{(j)}$ are appropriate to a Runge-Kutta solver of order P; in particular $0<b^{(j)}<1$, and $\,\sum_{j=1}^P b^{(j)} = 1$. We must estimate $y^{n+1}$ on an ad hoc basis; in this scheme we choose
\begin{align}
\label{macroy}y^{n+1} = y^{n,1}_{M_1} \;\;.
\end{align}

\noindent
The three numerical methods described above can now be considered as special cases of the HMM formulation \eqref{kj}--\eqref{macroy}, {and are illustrated in Figure~\ref{fig.expl}.}

\vspace{2pt}
\begin{enumerate}
\item \textbf{BA}: Choose $M_1 = 1$ and all other $M_j=0$. As described above, the BA takes smaller time steps in the slow variable than the HMM methods; we choose the macro step to be $\tDt$. We choose the number of iterations of the method $\tilde{n} $, so that the time elapsed over an application of the BA is $t^{\tilde{n}} = \tilde{n} \tDt$. The BA is illustrated in Figure~\ref{fig.BA}. 
\vspace{2pt}
\item \textbf{HMM1}: Fix $M_j = M > 0$ for all j, so that the micro solver is employed before each increment is estimated during the application of the macro solver. Choose the macro step $\Dt = M\tDt$ and number of iterations $n=\tilde{n}/M$, so that the time elapsed over an application of HMM1 is $t^{{n}} = {n} \Dt = \tilde{n}\tDt = t^{\tilde{n}}$ in the BA. The formulation HMM1 is illustrated in Figure~\ref{fig.HMM1}. 
\vspace{2pt}
\item \textbf{HMM2}: Choose $M_1 = M$ as above, and all other $M_j=0$, so that the micro solver is employed before the first increment is estimated in the macro solver, and held fixed during the macro solver. As in HMM1, choose macro step $\Dt$ and $n$ iterations so that $t^n = t^{\tilde{n}}$. The formulation HMM2 is illustrated in Figure~\ref{fig.HMM2}.
\end{enumerate}

\renewcommand{\thesubfigure}{\alph{subfigure}}

\begin{figure}
              \centering      
              \begin{subfigure}[h]{.5\textwidth}
                \centering
\includegraphics[width=\textwidth]{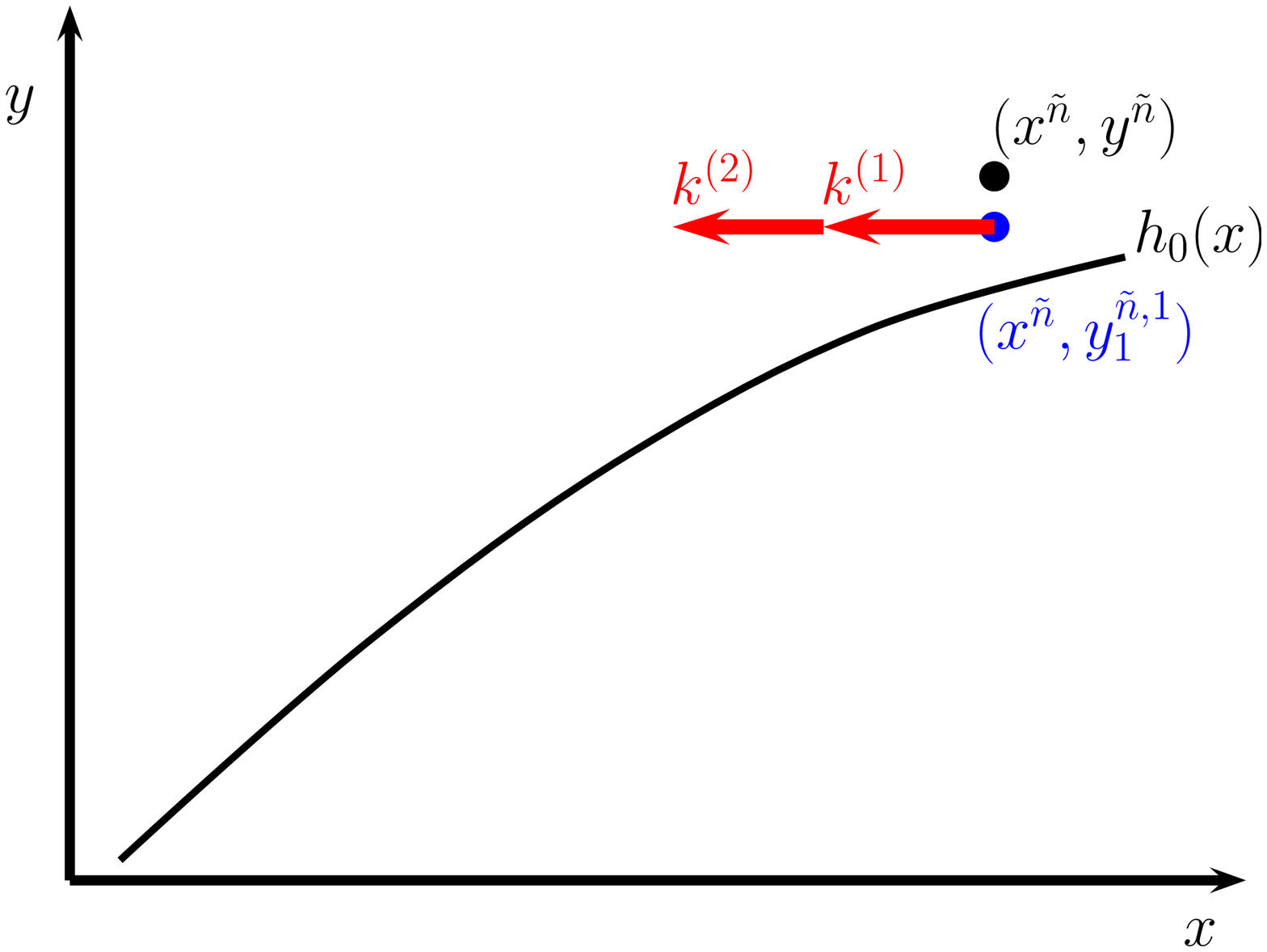}
\vspace{-20pt}
\caption{BA }
\label{fig.BA}
        \end{subfigure}
             
        \begin{subfigure}[h]{.5\textwidth}
                \centering
\includegraphics[width=\textwidth]{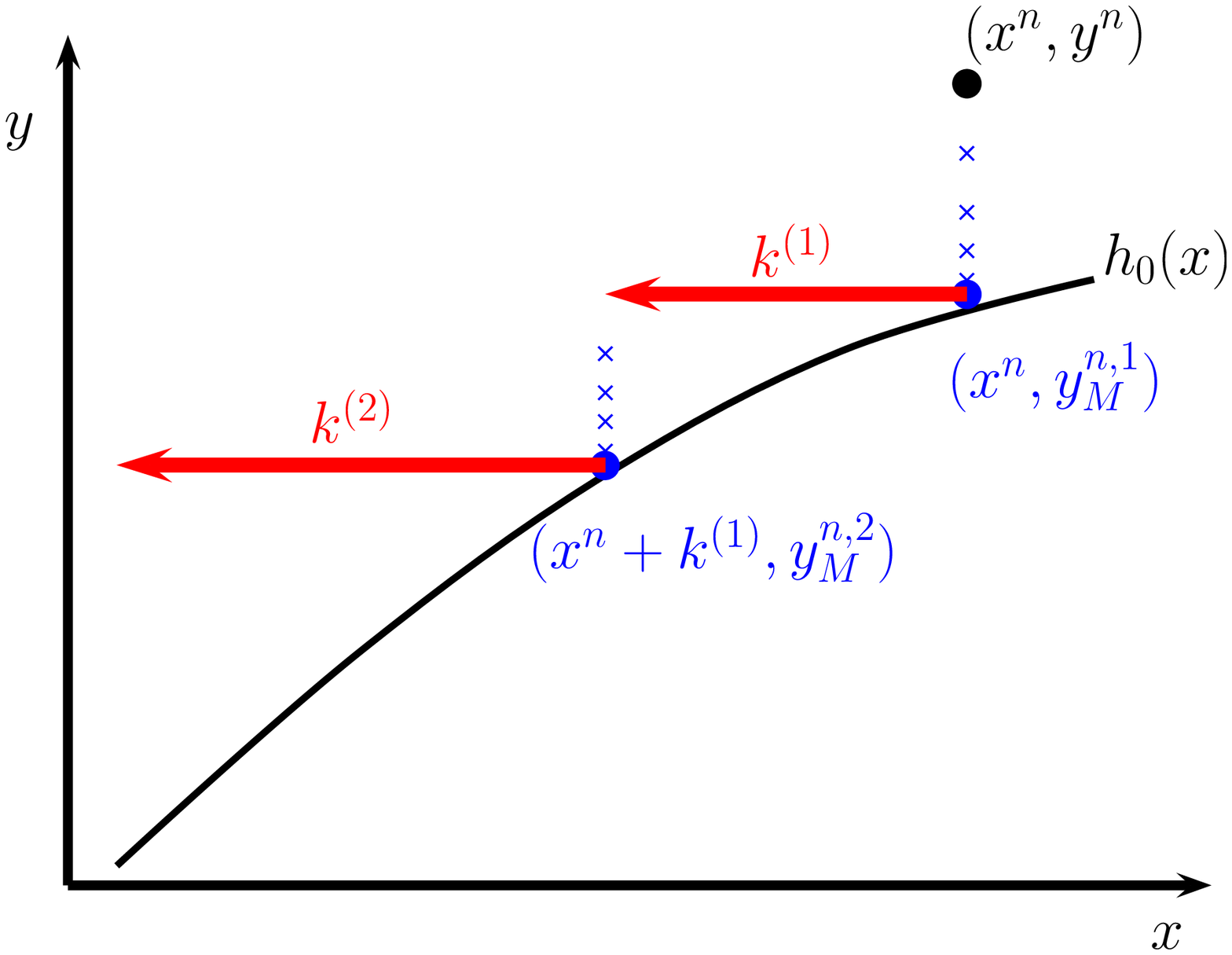}
\vspace{-20pt}
\caption{HMM1 }
\label{fig.HMM1}
        \end{subfigure}%
        \begin{subfigure}[h]{.5\textwidth}
                \centering
\includegraphics[width=\textwidth]{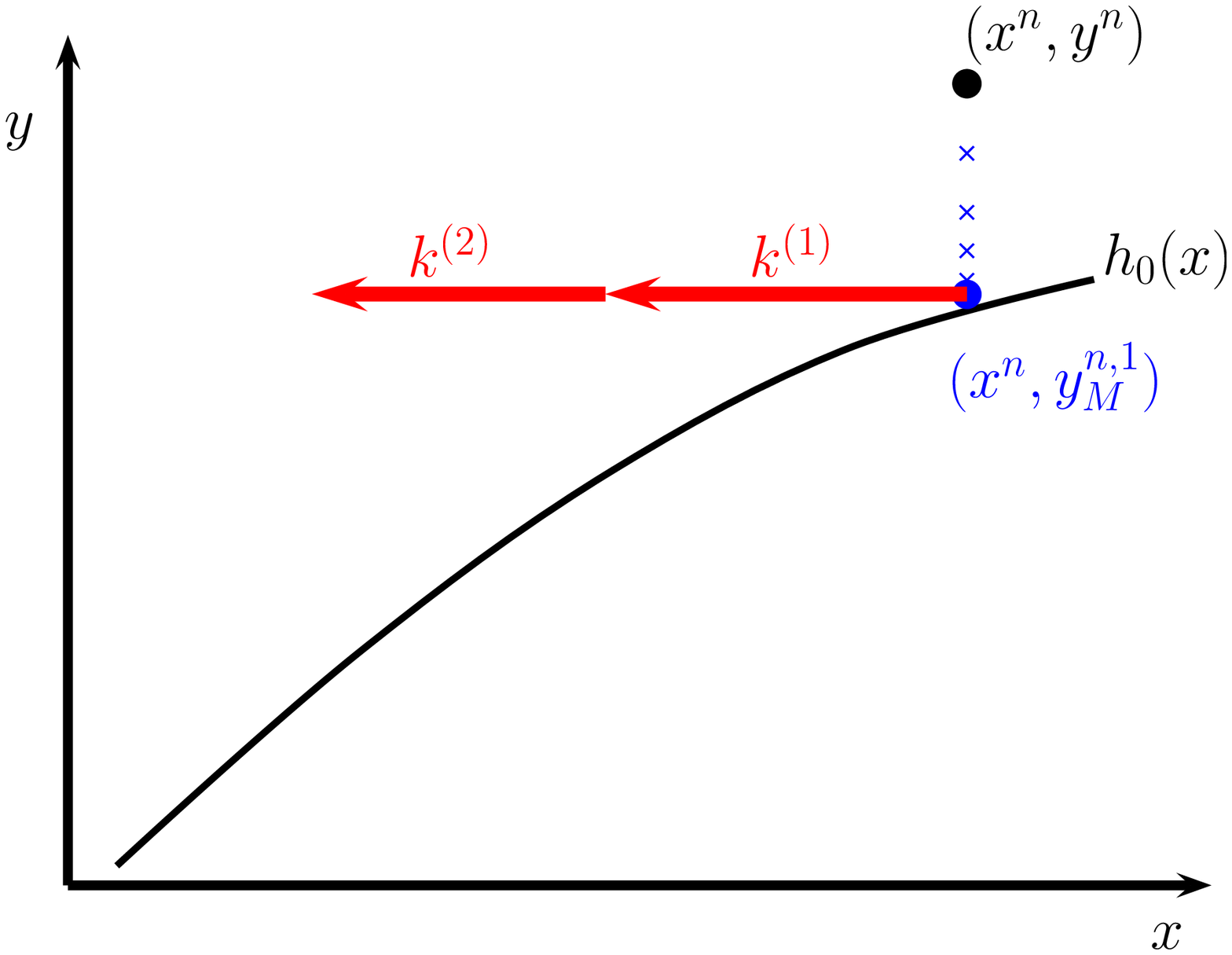}
\vspace{-20pt}
\caption{HMM2 }
\label{fig.HMM2}
       \end{subfigure}
        \caption{Sketch of the numerical methods for a second order Runge-Kutta macro solver.  For this scheme $a^{(2)} = 1$ and $b^{(1)}=b^{(2)}={1}/{2}$. Initial conditions are $(x^n,y^n)$ at the top right; applications of the micro solver are given by crosses while the increments are displayed as vectors. To focus on the different ways the increments $k^{(j)}$ are defined in the three formulations, we omit the final step of the macro solver, which is the weighted sum $x^{n+1} = x^n + (k^{(1)} + k^{(2)})/2$. }
                 \label{fig.expl}
\end{figure}

\section{Error Analysis}\label{sec:analysis}
We will establish bounds on the error between the slow variables $x^n$ given by the general formulation of HMM \eqref{kj}--\eqref{macroy} and the true value of the reduced slow dynamics $X(t^n)$. This provides error bounds for the Boosting Algorithm that improve on the bounds obtained in \cite{EEtAl09}, and furthermore establishes the distinction between the HMM1 and HMM2 formulations.\\

\noindent
We assume that the slow dynamics \eqref{baseslow} is Lipshitz continuous in the local region in which the slow manifold exists, so that there exist constants $L_f$ and $C_f$ satisfying
\begin{align*}
|f(x_1,y_1)-f(x_2,y_2)| \leq& L_f(|x_1-x_2|+|y_1-y_2|)\; , \\
|f(x,y)| \le& C_f \; .
\end{align*}
 Furthermore we assume that the approximate slow manifold $h_0(x)$ is Lipshitz continuous, so that there exists a constant $L_h$ satisfying 
 \begin{align*}
|h_0(x_1)-h_0(x_2)| \leq L_h|x_1-x_2|\; .
\end{align*}
Combining the previous two assumptions, the reduced dynamics \eqref{e.G0} is Lipshitz continuous, so that there exists a constant $L_F$ satisfying
 \begin{align*}
|F(X_1)-F(X_2)| \leq L_F|X_1-X_2|\; ,
\end{align*}
with $L_F \le L_f(1+L_h)$. Under these assumptions, we formulate the following Theorem.\\

\begin{theorem}
\label{theorem.main}
For the general formulation of HMM described by \eqref{kj}--\eqref{macroy}, with Runge-Kutta macro solver of order $P>1$ and micro solver of order $p\ge 1$, and for any $T\ge 0$, there exists a constant $C>0$ such that the error $|x^n - X(t^n)|$ is bounded for $t \le T$ by
\begin{align*}
|x^n - X(t^n)| \le C \left( \Dt^P + \rho^{\min(M_j)}\!\!\left(-\frac{\dt}{\eps}\right) \left| d^{n,\max}\right| + \eps \right) \; ,
\end{align*}
where $\rho{(-\dt/\eps)}$ is the linear amplification factor for the micro solver {and is derived in the proof on p. 9,} and $|d^{n,\max}|$ measures the maximal deviation of the fast variable from the slow manifold over all macro steps. Here $0\le \rho(-\dt/\eps) < 1$, and $\rho^{\min(M_j)}(-\dt/\eps) {= [\rho (-\delta t/\epsilon)]^{\min(M_j)}}$ signifies the exponential attraction of the fast dynamics towards the slow manifold over the microsteps. 
\end{theorem}

\medskip
\noindent
In order to apply this result to the formulations BA, HMM1 and HMM2, we employ the following practical assumptions.\\
We assume the parameters $\dt,\, M,\, \Dt$ of the HMM and BA solvers satisfy
\begin{align*}
\dt \le \eps \le M\dt < \tDt < \Dt <1 \;\;.
\end{align*}
We assume the micro solver accurately resolves the fast variable compared to the drift in the slow manifold over an increment or macro step, with
\begin{align}
 \label{pa.HMM}\rho^{M}\!\!\left(-\frac{\dt}{\eps}\right)|y^n - h_0(x^n)| <& L_h C_f \Dt  \;\; 
\end{align}
in HMM1 and HMM2, and 
\begin{align}
 \label{pa.BA}\rho^{\tilde{n}}\!\!\left(-\frac{\dt}{\eps}\right)\left|y^0 - h_0(x^0)\right| <&  \frac{L_h C_f \tDt \, \eps}{\dt} \;\;
\end{align}
 in the BA.\\ 

\noindent
Theorem~\ref{theorem.main} then has the following implications for the three specific implementations of HMM described in Section~\ref{sec:methods}:

\begin{corollary}
\label{c.BA}
For the Boosting Algorithm,
\begin{align*}
|x^{\tilde{n}} - X(t^{\tilde{n}})| & \le C\left(\tDt^P + \tDt\frac{\eps}{\dt} + \eps \right) \;.
\end{align*}
In particular, the predicted error is $|x^n - X(t^n)| \sim \tDt$ regardless of the order of the macro solver.\\
\end{corollary}
\noindent
 {The interpretation of the BA as a rescaling of time, described in Section \ref{sec:methods}, was used in \cite{EEtAl09} to derive an error bound for the BA. In this interpretation, the BA is thought of as a numerical scheme applied to \eqref{baseslow}--\eqref{basefast} with a larger value of $\eps$. Consequently the error incurred by the BA "consists of two parts: The error due to boosting the parameter $\eps$ and the error due to numerical solution of the boosted model" (\cite{EEtAl09}, p.5444). This bound can be written as
\begin{align*}
|x^{\tilde{n}} - X(t^{\tilde{n}})| & \le C\left( \tDt{\eps}/{\dt} + \tDt^P  + \left(\dt/\eps\right)^p \right) \;.
\end{align*}
That is, the error bound presented in \cite{EEtAl09} predicts the linear scaling of the error with $\tDt$, but also} contains an additional term proportional to $(\dt/\eps)^p$. This term would restrict the microstep to $\dt \ll \eps$ and suggest that a higher order micro solver might be employed with $p>1$, {otherwise the BA would incur $\mathcal{O}(1)$ error stemming from the micro solver}. Our bound does not have this limitation. In particular, one can choose larger $\dt$ to rapidly propagate the fast variable towards the slow manifold, and a first order micro solver with $p=1$ to lessen the computational cost, without increasing the predicted error.\\

\begin{corollary}
\label{c.HMM1}
For the Heterogeneous Multiscale Methods in the formulation HMM1,
\begin{align*}
|x^n - X(t^n)| \le C \left( \Dt^P + \Dt\, \rho^{M}\left(-\frac{\dt}{\eps}\right) + \eps \right) \; .
\end{align*}
\end{corollary}
\medskip
\begin{corollary}
\label{c.HMM2}
For the Heterogeneous Multiscale Methods in the formulation HMM2,
\begin{align*}
|x^n - X(t^n)| &\le C\left(  \Dt^P + \Dt+ \eps \right)\; .
\end{align*}
\end{corollary}
\noindent
In Corollary~\ref{c.HMM1}, the error bound may be dominated by any of the three terms. However, if $\eps \ll \Dt^P$ and $\Dt \, \rho^M(-\dt/\eps) \ll \Dt^P$, the error bound for the HMM1 formulation is identical to the bound for a Runge-Kutta solver applied to the reduced system. In Corollary~\ref{c.HMM2}, employing the practical assumption $\eps < \Dt <1$, the error bound is dominated by the linear term in $\Dt$ regardless of the order of the macro solver. \\
Some descriptions of HMM appear to recommend the HMM2 formulation by describing the application of the micro and macro solvers separately, as discussed earlier. However, Corollaries~\ref{c.HMM1} and~\ref{c.HMM2} outline the clear advantages of HMM1 over HMM2 for higher-order macro solvers.\\

\noindent
We now prove Theorem~\ref{theorem.main}.
\begin{proof}
We follow the procedure established in \cite{EEngquist03} for bounding the error in HMM methods, and in particular we make use of the convergence proof for our system \eqref{baseslow}--\eqref{basefast} in \cite{E03}. \\
Let us consider the numerical approximation of $X(t^n)$, the solution to the reduced slow system \eqref{e.G0} at time $t^n$, by a single step of a Runge-Kutta method of order P with initial condition $X(t^{n-1})$. The increments are given by
\begin{align}
\label{Kj} {K}^{(j)}\left((X(t^{n-1})\right)& = \Dt \, F \left(X(t^{n-1}) + a^{(j)} K^{(j-1)} \right) \;\; ,
\end{align}
and the solver is then given by the weighted sum
\begin{align}
\label{macroX} X^{n} &= X(t^{n-1}) + \sum_{j=1}^P b^{(j)} {{K}^{(j)}\left(X(t^{n-1})\right)}{}\;\; .
\end{align}
The Runge-Kutta solver approximates the reduced dynamics to order $P+1$ over a single step (see for instance \cite{Iserles}), so that 
\begin{align}
\label{RK.step}X(t^n) = X(t^{n-1}) + \Dt \hat{F}\left(X(t^{n-1})\right)+ \mathcal{O}(\Dt^{P+1}) \;,
\end{align}
where we write the vector field of the solver as
\begin{align}
\label{Fhat} \hat{F}\left(X(t^{n-1})\right) &= \sum_{j=1}^P b^{(j)} \frac{{K}^{(j)}\left(X(t^{n-1})\right)}{\Dt} \;\; .
\end{align}
Similarly for the HMM formulation \eqref{kj}--\eqref{macroy} we write the vector field of the macro solver as
\begin{align}
\label{fhat} \hat{f}(x^n) &= \sum_{j=1}^P b^{(j)} \frac{{k}^{(j)}(x^{n},y^n)}{\Dt} \;\; ,
\end{align}
so that
\begin{align}
\label{macrof} x^{n+1} &= x^n + \Dt \hat{f}(x^n) \; .
\end{align}

\noindent
We expand $x^n - X(t^n)$ by substituting \eqref{RK.step} and \eqref{macrof}, obtaining
\begin{align*}
x^n - X(t^n) =& x^{n-1} - X(t^{n-1})+ \Dt \left(\hat{f}(x^{n-1},y^{n-1})  - \hat{F}\left(X(t^{n-1})\right)\right) + \mathcal{O}(\Dt^{P+1})\\
=& x^{n-1} - X(t^{n-1}) +\Dt \left(\hat{F}(x^{n-1}) - \hat{F}\left(X(t^{n-1})\right)\right)\\
&+ \Dt \left(\hat{f}(x^{n-1},y^{n-1}) - \hat{F}(x^{n-1})\right) + \mathcal{O}(\Dt^{P+1}) \\
=& \left(x^{n-1} - X(t^{n-1})\right)\left[1 +\Dt \hat{F}\left(X(t^{n-1}) + \theta (x^{n-1} - X(t^{n-1}))\right) \right]\\
&+ \Dt \left(\hat{f}(x^{n-1},y^{n-1}) - \hat{F}(x^{n-1})\right) + \mathcal{O}(\Dt^{P+1}) \;\;,
\end{align*}
where we have applied the Mean Value Theorem with $0\le \theta \le 1$. Applying absolute values and employing the Lipshitz continuity of the reduced dynamics, we obtain 
\begin{align*}
\left| x^n - X(t^n)\right|\le & \left(1 + L_F \Dt \right)\left|x^{n-1} - X(t^{n-1})\right| + \Dt \left|\hat{f}(x^{n-1},y^{n-1}) - \hat{F}(x^{n-1})\right| \\
&+C\Dt^{P+1}\;\;.
\end{align*}
Iterating with initial condition $X(0) = x^0 $ and $t^n \le T$ yields
\begin{align}
\label{xpause} |x^n - X(t^n)| &\le \frac{e^{L_F T}}{L_F} \left( C \Dt^P + \max_{1\le i < n} \left| \hat{f}(x^i,y^i) - \hat{F}(x^i)\right|  \right) \;\; .
\end{align}
This is the standard bound for HMM methods, established for very general systems in \cite{EEngquist03}. The first term is the error associated with a Runge-Kutta solver applied to the reduced system. The second term measures how well the HMM method approximates this Runge-Kutta solver, and is labelled $e(\text{HMM})$. 
\\
We now bound $|\hat{f}(x^n,y^n) - \hat{F}(x^n)| $. Substituing \eqref{Fhat} and \eqref{fhat} we obtain
\begin{align}
\nonumber \left| \hat{f}(x^n,y^n) - \hat{F}(x^n)\right| &\le  \sum_{j=1}^P \frac{b^{(j)}}{\Dt}  \left| {{k}^{(j)}(x^{n},y^n)-  {K}^{(j)}(x^{n})}{} \right| \\
\label{fpause} &\le \frac{1}{\Dt}  \max_{1\le j\le P}\left| {{k}^{(j)}(x^{n},y^n)-  {K}^{(j)}(x^{n})}{}\right| \;\;,
\end{align}
on employing the weighting condition $\sum_{j=1}^P b^{(j)} =1$.\\ We bound $\left| {{k}^{(j)}(x^{n},y^n)-  {K}^{(j)}(x^{n})}{}\right|$, substituting \eqref{kj} and \eqref{Kj} for the increments with
\begin{align}
\label{uncomf} \left| {{k}^{(j)}(x^{n},y^n)}{} -  {{K}^{(j)}(x^{n})}{} \right|=&  \Dt \left| f(x^{n}+ a^{(j)} k^{(j-1)},y^{n,j}_{M_j})-  F(x^{n}+ a^{(j)} K^{(j-1)}) \right| \\
\nonumber =&  \Dt \big| f(x^{n}+ a^{(j)} k^{(j-1)},y^{n,j}_{M_j})\\
\nonumber &\quad\,\,-  f(x^{n}+ a^{(j)} K^{(j-1)},h_0(x^{n}+ a^{(j)} K^{(j-1)})) + \mathcal{O}(\eps) \big| \;\; ,
\end{align}
where we used the zeroth order approximation of the slow reduced vector field \eqref{e.G0}. Employing the Lipshitz continuity of the slow dynamics and slow manifold yields
\begin{align}
\nonumber \left| {{k}^{(j)}(x^{n},y^n)}{} -  {{K}^{(j)}(x^{n})}{} \right| \le&  L_f \Dt \Big(\left|a^{(j)} k^{(j-1)} - a^{(j)} K^{(j-1)}\right| \\
\nonumber&\qquad\quad+ \left| y^{n,j}_{M_j} - h_0(x^{n}+ a^{(j)} K^{(j-1)}) \right|  \Big)  + C \Dt\eps \\
\nonumber \le&  L_f \Dt \Big(a^{(j)}\!\!\left| k^{(j-1)} -  K^{(j-1)}\right|\! + \!\left| y^{n,j}_{M_j} - h_0(x^{n}+ a^{(j)} k^{(j-1)}) \right|   \\
\nonumber &+ \left|h_0(x^{n}+ a^{(j)} k^{(j-1)})  - h_0(x^{n}+ a^{(j)} K^{(j-1)}) \right|  \Big) + C \Dt\eps \\
\nonumber \le&  L_f (1+ L_h) \Dt a^{(j)}\left|k^{(j-1)} -  K^{(j-1)}\right| \\
\label{kpause} &+ L_f \Dt\left| y^{n,j}_{M_j} - h_0(x^{n}+ a^{(j)} k^{(j-1)}) \right|   +  C \Dt\eps \;.
\end{align}
Defining the distance between the fast variable $y^{n,j}_m$ and the approximate slow manifold $h_0(x^n+a^{(j)}k^{(j-1)})$ by
\begin{align*}
d^{n,j}_{m} = y^{n,j}_{m} - h_0(x^{n}+ a^{(j)} k^{(j-1)}) \;\; ,
\end{align*}
we write \eqref{kpause} with $j=1$ as
\begin{align}
\label{hmf} \left| {{k}^{(1)}(x^{n},y^n)}{} -  {{K}^{(1)}(x^{n})}{} \right| \le
&  L_f \Dt | d^{n,1}_{M_1}  |   +  C \Dt\eps \;.
\end{align}
Iterating \eqref{kpause}, seeded with \eqref{hmf}, we obtain
\begin{align*}
 \left| {{k}^{(j)}(x^{n},y^n)}{} -  {{K}^{(j)}(x^{n})}{} \right| \le& L_f \Dt \max_{1\le j \le P} | d^{n,j}_{M_j} | + C \Dt\eps+ \mathcal{O}(\Dt^2\max_{1\le j \le P}| d^{n,j}_{M_j} |,\Dt^2 \eps)    \;.
\end{align*}
On substitution into \eqref{fpause} we obtain
\begin{align}
\label{fpause2} \left| \hat{f}(x^n,y^n) - \hat{F}(x^n)\right|  &\le L_f  \max_{1\le j \le P}| d^{n,j}_{M_j} | + C \eps + \mathcal{O}(\Dt\max_{1\le j \le P}| d^{n,j}_{M_j} |)   \;\;.
\end{align}
We have expressed our bound over the vector fields $\hat{f}$ and $\hat{F}$ as a bound over the micro steps. The terms $L_f \max_{1\le j \le P}| d^{n,j}_{M_j} | + C\eps $ measure via \eqref{uncomf} the mismatch between the slow vector field $f$ after an applicaton of the micro solver and the reduced vector field $F$. We now bound $|d^{n,j}_{M_j} |$.\\

\noindent
Recall that $y^{n,j}_{M_j}$ is produced by the micro solver, which resolves the fast dynamics \eqref{basefast} with the slow variable held fixed. {At any fixed $x$, define $d_\eps = y_\eps - h_0(x)$}, with linear dynamics given by 
\begin{align}
\nonumber \dot d_\eps &= \dot y_\eps\\
\label{d dot} & = - \frac{1}{\eps}d_\eps \;\; ,
\end{align}
{so that
\begin{align}
\label{d.lin}\frac{d^r \left( d_\eps \right)}{dt^r} &= \left(-\frac{1}{\eps}\right)^r d_\eps \;\;.
\end{align}
}
By construction, $d^{n,j}_{M_j} $ is the output of the micro solver applied to \eqref{d dot} with initial condition $d_\eps(0) = y^{n,j}_0 - h_0(x^n+a^{(j)}k^{(j-1)})$. {We now use this interpretation of $d^{n,j}_m$ to establish the rate at which various micro solvers propagate the fast variable to the slow manifold.\\
It is well known that a single step in a Runge-Kutta solver of order p, when Taylor expanded term by term, matches the Taylor series around the previous step to order p in $\dt$. That is, 
\begin{align*}
d^{n,j}_{m+1} &= d^{n,j}_{m} + \dt \,\left.\dot d_\eps\right|_{d^{n,j}_{m}} + \frac{\dt^2}{2!} \left.\ddot d_\eps\right|_{d^{n,j}_{m}} + \dots + \frac{\dt^p}{p!} \left.\frac{d^p (d_\eps)}{dt^p}\right|_{d^{n,j}_{m}} +\O(\dt^{p+1}) \;\;.
\end{align*}
With the linear dynamics of \eqref{d.lin}, the Taylor expansion is unnecessary: every increment in the Runge-Kutta step from $d^{n,j}_{m}$ to $d^{n,j}_{m+1}$ can be rewritten in terms of $d^{n,j}_{m}$ directly. Consequently there is no error term in the above equation, and
\begin{align}
 \nonumber d^{n,j}_{m+1} &= d^{n,j}_{m} + \dt \left.\dot d_\eps\right|_{d^{n,j}_{m}} + \frac{\dt^2}{2!} \left.\ddot d_\eps\right|_{d^{n,j}_{m}} + \dots + \frac{\dt^p}{p!} \left.\frac{d^p (d_\eps)}{dt^p}\right|_{d^{n,j}_{m}} \;\;.
 \end{align}
Employing \eqref{d.lin} we obtain the bound on $d^{n,j}_m$ as stated in \cite{E03},
\begin{align}
d^{n,j}_{m+1}&=\label{dpause}\rho\left(-\frac{ \dt}{\eps}\right) d^{n,j}_{m} \;\;,
\end{align}}
where
\begin{align*}
 \rho\left(-\frac{ \dt}{\eps}\right) & = \sum_{j=0}^p \frac{(-\frac{ \dt}{\eps})^j}{j!} \;\; .
\end{align*}
For the fast variable to converge to the (approximate) slow manifold, we require $|\rho(-\dt/\eps)| < 1$, which is assured for any $0 < \dt \le \eps$. Note that that from \eqref{dpause} we have $\lim_{m\rightarrow\infty}d^{n,j}_m = 0$, implying that the fast variable converges to the approximate slow manifold $h_0$ over the micro solver, rather than converging to the true slow manifold $h_\eps$. Consequently, \eqref{fpause2} assures that the HMM vector field $\hat{f}$ converges to the zeroth order approximation of the reduced dynamics. \\
Substituting \eqref{dpause} into \eqref{fpause2}, and \eqref{fpause2} into 
\eqref{xpause}, we obtain the main result of Theorem~\ref{theorem.main}, that
\begin{align}
\label{xpause2} |x^n - X(t^n)| \le C \left( \Dt^P + \rho^{\min(M_j)}\!\!\left(-\frac{\dt}{\eps}\right) \left| d^{n,\max}\right| + \eps \right) \;\; ,
\end{align}
where we defined the maximal deviation $\left| d^{n,\max}\right|$ of the fast variable from the slow manifold over the applications of the macro solver,
\begin{align*}
\left| d^{n,\max}\right| =\max_{\substack{0 \le i <n \\ 1 \le j \le P}} | d^{i,j}_{0} | \;\;.
\end{align*}
Equation \eqref{xpause2} already contains the term that differentiates HMM1 from HMM2: in HMM1, $\min(M_j) = M$, while in HMM2 (and in BA), $\min(M_j) = 0$.
\end{proof}

\medskip
\noindent
We now prove Corollary~\ref{c.BA}, \ref{c.HMM1} and \ref{c.HMM2} by bounding $| d^{n,\max}|$ in each of the three formulations. 

\begin{proof}
First fix $n$ and consider $|d^{n,j}_0|$ over the increments. For $j>1$ we employ the Lipshitz continuity of the slow dynamics to bound 
\begin{align}
\nonumber |d^{n,j}_0| =& | y^{n,1}_{M_1} - h_0(x^n + a^{(j)} k^{(j-1)}) | \\
\nonumber \le& | h_0(x^n) - h_0(x^n + a^{(j)} k^{(j-1)}) | + |y^{n,1}_{M_1} - h_0(x^n)| \\
\label{d.inc} \le & L_h C_f \Dt + \rho^{M_1}\left(-\frac{\dt}{\eps}\right) |d^n| \;\;,
\end{align}
where we have employed $a^{(j)}\le 1$, and where we have defined the distance of the fast variable from the slow manifold at the beginning of the $n$-th macro step,
\begin{align*}
 |d^n| = |y^n - h_0(x^n)| \;.
 \end{align*} 
For $j=1$, $|d^{n,1}_0| = |d^n| $. Analogously to \eqref{d.inc} we now bound $|d^n|$, employing \eqref{macrox} and \eqref{macroy} to obtain
\begin{align}
\nonumber |d^{n+1}| =& | y^{n,1}_{M_1} - h_0(x^{n+1} ) | \\
\nonumber \le& | h_0(x^n) - h_0(x^{n+1}) | + |y^{n,1}_{M_1} - h_0(x^n)| \\
\label{d.mac} \le & L_h C_f \Dt + \rho^{M_1}\left(-\frac{\dt}{\eps}\right) |d^n| \;\;.
\end{align}
 We note that the bounds for the increments \eqref{d.inc} and the macro steps \eqref{d.mac} are the same. In each bound, the first term measures the changing value of the slow manifold over the drift of the slow variable, and the second term measures the distance of the fast variable from the approximate slow manifold after the last macro step. \\
We will now use \eqref{d.inc} and \eqref{d.mac} to bound $|d^{n,\max}|$ for the three formulations HMM1, HMM2 and BA. \\

\noindent
Recall that for HMM1 and HMM2, $M_1 = M$. For these formulations, we employ the practical assumption \eqref{pa.HMM} that the micro solver accurately resolves the fast variable compared to the drift in the slow manifold over an increment or macro step, which we recall as
\begin{align*}
\rho^{M}\left(-\frac{\dt}{\eps}\right)|d^n| <& L_h C_f \Dt \;\; .
\end{align*}
Then we obtain to lowest order the bound for HMM1 and HMM2, 
\begin{align}
\label{dnmax1} \left| d^{n,\max}\right| \le 2L_h C_f \Dt \;\; .
\end{align}
For the BA, with $M_1=1$, we cannot neglect $\rho^{1}\left(-\frac{\dt}{\eps}\right) |d^{\tilde{n}}|$; instead we iterate \eqref{d.mac} to obtain
\begin{align*}
\left|d^{\tilde{n}}\right| 
\le &\rho^{\tilde{n}}\left(-\frac{\dt}{\eps}\right)\left|d^0\right| + L_h C_f \tDt\frac{1-\rho^{\tilde{n}}\left(-\frac{\dt}{\eps}\right)}{1-\rho\left(-\frac{\dt}{\eps}\right)} \\
 \le&\rho^{\tilde{n}}\left(-\frac{\dt}{\eps}\right)\left|d^0\right| +  L_h C_f \tDt\frac{1}{1-\rho\left(-\frac{\dt}{\eps}\right)} \\
 \le&\rho^{\tilde{n}}\left(-\frac{\dt}{\eps}\right)\left|d^0\right| +  \frac{L_h C_f \tDt \, \eps}{\dt} \;\; .
\end{align*}
Recalling that $\tilde{n} = nM$, we can now employ the practical assumption \eqref{pa.BA} that the applications of the micro solver dampen any initial distance of the fast variable from the slow manifold compared to the drift in the slow manifold over an increment or macro step, which we recall as
\begin{align*}
\rho^{nM}\left(-\frac{\dt}{\eps}\right)\left|d^0\right| <&  \frac{L_h C_f \tDt \, \eps}{\dt} \;\;,
\end{align*}
yielding to lowest order the bound
\begin{align}
 \left| d^{n,\max}\right| 
\label{dnmax2} \le&  2\frac{L_h C_f \tDt \, \eps}{\dt} \;\; .\\
\nonumber
\end{align}

\noindent
Substituting \eqref{dnmax2} into \eqref{xpause2} with $\min(M_j) = 0$ establishes Corollary~\ref{c.BA}; similarly \eqref{dnmax1} with $\min(M_j) = M$ establishes Corollary~\ref{c.HMM1}, and \eqref{dnmax1} with $\min(M_j) = 0$ establishes Corollary~\ref{c.HMM2}.
\end{proof}

\section{Numerics}
\label{sec:numerics}
We now illustrate the key predictions of Corollaries \ref{c.BA}, \ref{c.HMM1}, and \ref{c.HMM2}, which we recall here with a second-order Runge-Kutta macro solver, i.e. $P=2$, and a forward Euler micro solver, i.e. $\rho(-\dt/\eps) = 1-\dt/\eps$. We present the result for each method along with the predicted behaviour of the error $|x^n - X(t^n)|$ as the parameters $\Dt,\, \eps,\, \dt,\, M$ are varied. 
\vspace{2pt}
\begin{enumerate}
\item BA: 
\begin{align*}
 |x^n - X(t^n)| \le&  C\left( \left(\frac{\Dt}{M}\right)^2 + \frac{\Dt\, \eps}{M\, \dt}  + \eps\right)   \;\;,
\end{align*}
where we have used $\tDt = \Dt/M$ so that the BA may be compared to HMM1 and HMM2. So long as the practical assumptions $M\dt < \Dt < 1$, $\dt \le \eps$ are satisfied, the error is always dominated by the second term; in particular the error is linear in $\Dt$ and in $\eps$. 
\vspace{2pt}
\item HMM1: 
\begin{align*}
 |x^n - X(t^n)| \le&  C\left(\Dt^2 + \Dt\left(1-\frac{\dt}{\eps}\right)^M+ \eps   \right) \;\;.
\end{align*}
Depending on the parameters of the method and system, any of the three terms could be the dominant error term. 
\vspace{2pt}
\item HMM2: 
\begin{align*}
 |x^n - X(t^n)| \le&  C\left(\Dt^2 + \Dt+ \eps   \right) \;\;.
\end{align*}
So long as the practical assumption $\eps < \Dt <1$ is satisfied, the error is always dominated by the second term; we expect the error to scale linearly with $\Dt$. Employing the practical assumption $\eps < M\dt$, we see that the second term in the bound for HMM2 is larger than the second term in the bound for the BA.
\end{enumerate}
\medskip

\noindent
To illustrate these results we present simulations of the Michaelis-Menten system employed in \cite{GearEtAl05},
\begin{align}
\label{mm_x}\dot{x}_{\eps} &= -x_{\eps} + (x_{\eps}+0.5)y_{\eps} \\
\label{mm_y}\dot{y}_{\eps} &= \frac{x_{\eps}- (x_{\eps}+1)y_{\eps}}{\eps} \; ,
\end{align}
with stable fixed point at $(0,0)$. We comment that the fast dynamics of this system is more complicated than our form \eqref{basefast}, but still converges to the slow manifold 
\begin{align}
\label{heps} h_\eps(x) &= \frac{x}{x+1}+\eps\frac{x}{2(x+1)^4} + \O(\eps^2) \;\;
\end{align}
exponentially quickly, so that $x_\eps$ converges to the slow reduced system
\begin{align}
\label{mm_X}\dot{X} &= -X + (X+0.5)h_\eps(X) \;\;.
\end{align}
We use initial conditions $x_{\eps}(0) = 1$, $y_{\eps}(0) = h_\eps(1)$. Under these conditions, the fast variable in the micro solver still converges to the slow manifold according to
\begin{align*}
\left|d^{n,j}_{m+1}\right| &\le \left(1-\frac{\dt}{\eps}\right) \left|d^{n,j}_{m}\right| \; ,
\end{align*}
so that Theorem~\ref{theorem.main} and the Corollaries apply. \\

\noindent
We present results for three situations. First we choose parameters so that the error in HMM1 is dominated by the term proportional to $\Dt^2$. We then demonstrate situations in which the error in HMM1 is dominated by the term proportional to $\Dt\left(1-{\dt}/{\eps}\right)^M$, and the term proportional to $\eps$.
\subsection{HMM1 with ideal parameters}
We consider the scaling of the error with the macro step $\Dt$ for parameters such that HMM1 is as accurate as a second-order Runge-Kutta solver applied to the reduced system, i.e. $|x^n - X(t^n)| \sim \Dt^2$. We choose $\eps=10^{-5}$ so that the error is not dominated by the term proportional to $\eps$. We choose $\dt = 0.2\eps$ and $M=30$ so that the fast variable converges very close to the approximate slow manifold with $(1-\dt/\eps)^M = 0.001  \ll \Dt$. The number of macro steps ranges from $10$ to $500$ to keep $T=5$ fixed. Figure~\ref{Dt1} clearly illustrates the identical quadratic scalings of the error in HMM1 and the error in a second-order Runge-Kutta solver applied to the reduced system, and the linear scaling of the error in HMM2 and the BA with $\Dt$. 

\begin{figure}
\includegraphics[width=\textwidth]{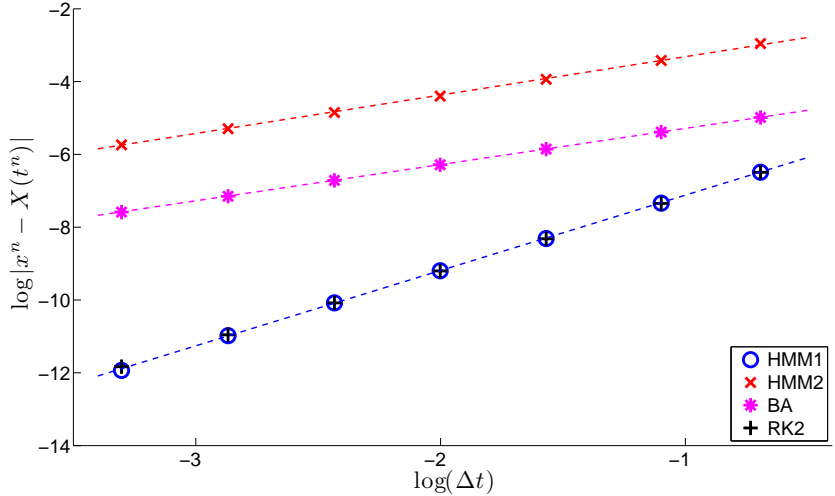}
\caption{Plot of the error versus $\log (\Dt)$ for fixed time of integration $t^n=5$ of the system \eqref{mm_x}-\eqref{mm_y}. The crosses labelled 'RK2' are from second-order Runge-Kutta simulations of the reduced system \eqref{mm_X}. The dashed lines are linear regressions with slopes of $2.07$ for HMM1, $1.06$ for HMM2 and $1.00$ for the BA.} 
\label{Dt1}
\end{figure}

\subsection{HMM error dominated by under-resolved fast variable}
We now consider the scaling of the error with the macro step $\Dt$ when the fast variable does not converge close to the slow manifold, so that the error in HMM1 is dominated by $\Dt (1-\dt/\eps)^M$; that is, the error scaling is linear in $\Dt$. We choose as before $\eps=10^{-5}$ and $\dt = 0.2\eps$, now with $M=10$ so that $(1-\dt/\eps)^M = 0.1  \approx \Dt$. The number of macro steps again ranges from $10$ to $500$ to keep $t^n=5$ fixed. We illustrate the linear error scaling with $\Dt$ of all three formulations in Figure~\ref{Dt2}.

\begin{figure}
\includegraphics[width=\textwidth]{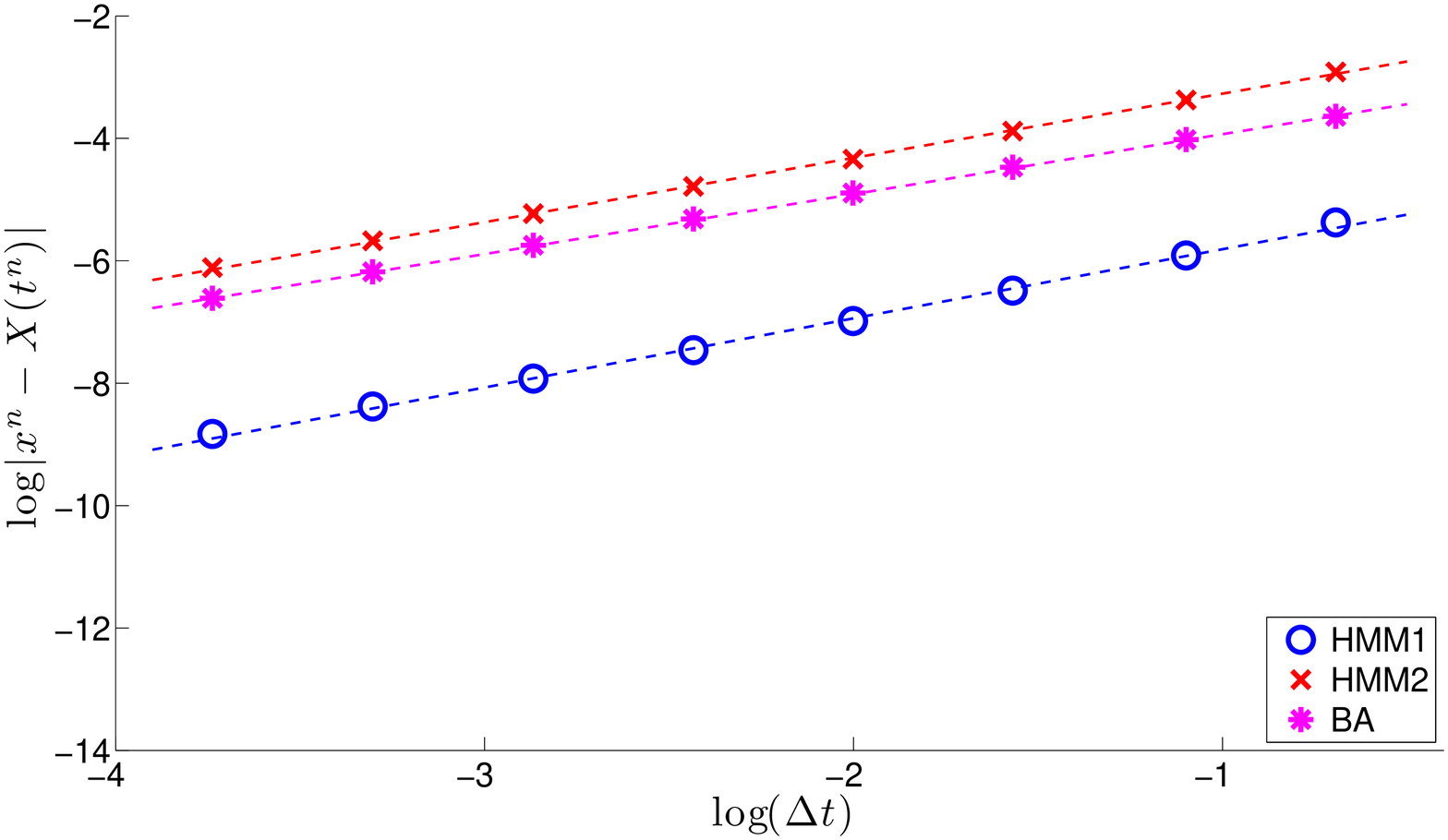}
\caption{Plot of the error versus $\log (\Dt)$ for fixed time of integration $t^n=5$ of the system \eqref{mm_x}-\eqref{mm_y}. The dashed lines are linear regressions with slopes of $1.13$ for HMM1, $1.05$ for HMM2 and $0.98$ for the BA.} 
\label{Dt2}
\end{figure}

\subsection{HMM error dominated by the scale separation}
We now consider the scaling of the error with $\eps$. This term is dominant in HMM1 when $\Dt^P \ll \eps$ and the fast variable converges well to the approximate slow manifold. In this situation, as discussed earlier, HMM1 converges to the zeroth order approximation of the reduced system $f(x,h_0(x))$ and the error is dominated by the $\O(\eps)$ difference between the approximate reduced system and the true reduced system $f(x,h_\eps(x))$. We choose as before $\dt = 0.2\eps$ and $M=30$ so that as before the fast variable converge very close to the approximate slow manifold. We use $\Dt = 0.1$ and take $n=50$ macro steps. We illustrate the linear error scaling with $\eps$ of the HMM1 and BA for $\eps > \Dt^2 = 0.01$ in Figure~\ref{Dt2}. For HMM2, the error term proportional to $\Dt$ always dominates the error so long as the practical assumption $\Dt \gg \eps$ is satisfied, as noted in Corollary~\ref{c.HMM2}. 
 
\begin{figure}
\includegraphics[width=\textwidth]{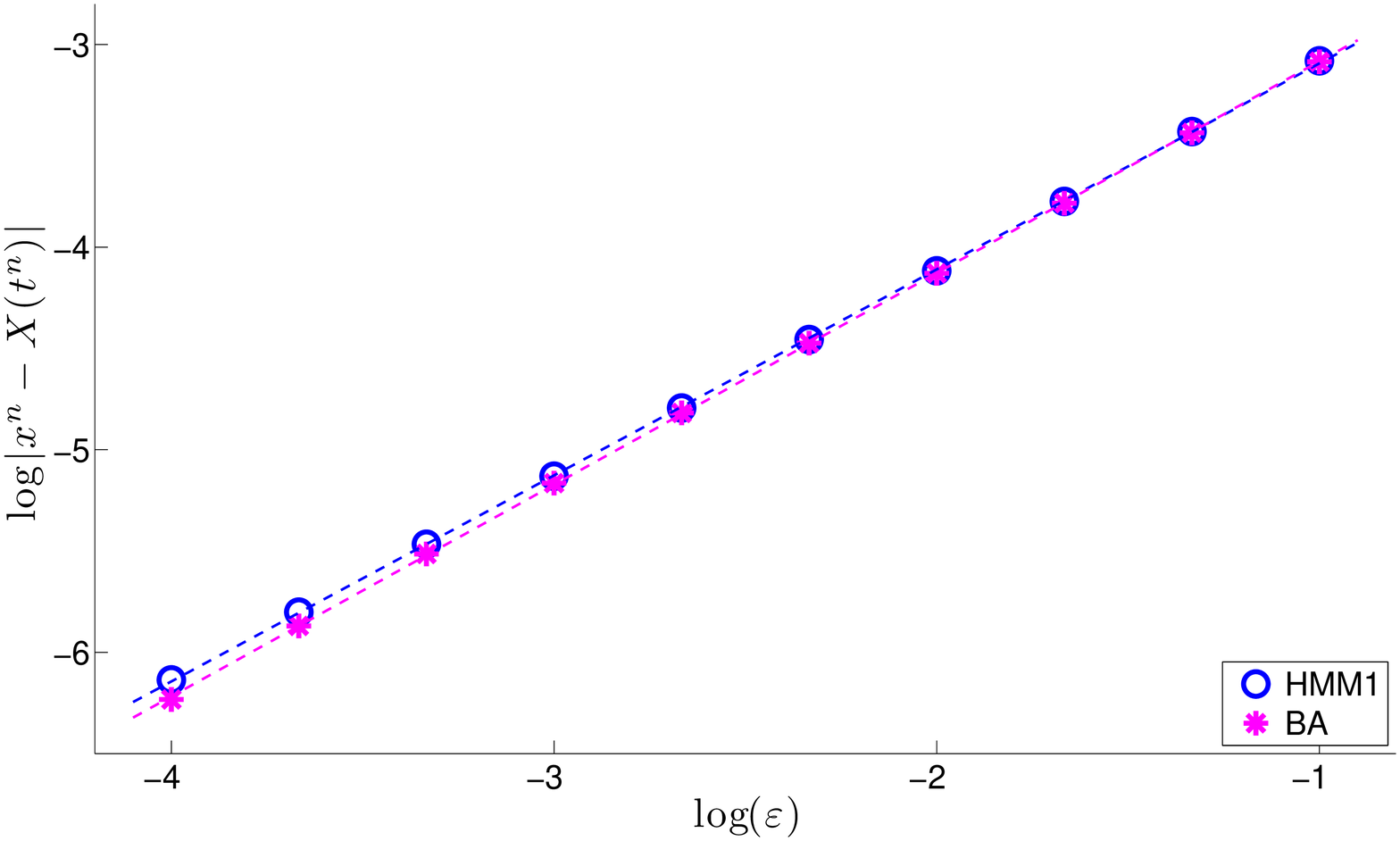}
\caption{Plot of the error versus $\log (\eps)$ for fixed time of integration $t^n=5$ of the system \eqref{mm_x}-\eqref{mm_y}. The dashed lines are linear regressions with slopes of $1.02$ for HMM1 and $1.05$ for the BA.} 
\label{eps}
\end{figure}

\section{Discussion} \label{sec:discussion}
\noindent
We have extended the proof in \cite{E03} to provide convergence results for the Heterogeneous Multiscale Methods with a Runge-Kutta macro solver and a variable number of micro steps. These results naturally lead to error bounds for three distinct numerical methods, all of which are special cases of the general formulation of HMM. In particular, we obtained an improved error bound for the Boosting Algorithm over the bound presented in \cite{EEtAl09}. The bound in \cite{EEtAl09} contains an additional term that suggests a restriction on the micro step to inefficiently small values or on the micro solver to high order solvers; the bound we derived does not have this restriction. \\

\noindent
We also obtained an error bound for two distinct formulations of HMM, and used these to demonstrate that the micro solver must be employed repeatedly during a single step of a Runge-Kutta macro solver, otherwise the error is dominated by a linear term in the macro step. If the micro solver is employed before every vector field evaluation in the macro solver as done in HMM1, the simulation of the slow variables is within $\O(\eps)$ of the output of a Runge-Kutta solver {of arbitrary order} applied to the reduced system. This provides a word of caution for the conceptual definitions of HMM. \\

\noindent
{While the HMM1 formulation is accurate, it is significantly more computationally expensive than HMM2 or the BA. Consequently, accurate formulations of HMM2 or the BA are still of interest. We now discuss some progress that has been made in designing such schemes. As we have shown in this paper, any such improvements will necessarily improve the estimate of the slow manifold used in the macro solver. \\
In \cite{LockerbyEtAl13}, a leapfrog coupling is shown to yield second order accuracy in the macro step for the HMM2 formulation. In particular, this allows the use of a second order implicit macro solver, which is difficult to implement in the HMM1 formulation and advantageous if the macro scale dynamics are also stiff. In the same paper a new numerical method is suggested, essentially a formulation of HMM2 with a varying number of micro steps before each macro time step. It is shown that this approach is advantageous in systems where the scale separation varies. \\
In \cite{LeeEngquist14} an extension of the BA is developed, called Variable Step size HMM. In this scheme the BA macro step varies smoothly from micro to macro values; it is shown that this approach enables the BA to achieve second order accuracy in the macro step under some conditions on the step size. \\
The bound we derived for the BA is partially applicable to VSHMM. In particular \cite{LeeEngquist14} uses fourth order Runge-Kutta for the micro solver in dissipative experiments when our analysis suggests one could gain in accuracy via faster convergence to the slow manifold, and save a factor of 4 computationally, by using forward Euler for the micro solver.}\\

\section*{Acknowledgements}
The author is grateful to Georg Gottwald for many useful comments. Any mistakes remain my own.

\bibliographystyle{siam}
\bibliography{bibliography}
\end{document}